\documentclass{amsart}

\usepackage[english]{babel}

\usepackage[letterpaper,top=2cm,bottom=2cm,left=3cm,right=3cm,marginparwidth=1.75cm]{geometry}

\usepackage{amsmath}
\usepackage{graphicx}
\usepackage{amsthm}
\usepackage{amsfonts}
\usepackage{comment}
\newtheorem{theorem}{Theorem}[section]
\newtheorem{corollary}[theorem]{Corollary}
\theoremstyle{definition}

\newtheorem{lemma}{Lemma}[section]
\usepackage{setspace}
\usepackage[colorlinks=true, allcolors=blue]{hyperref}
\date{}
\title[]{Recurrence relations for the coefficients of the confluent and Gauss hypergeometric functions in the complex plane}
\author{Zi-Qiao Xu, Zhong-Xuan Mao, Jing-Feng Tian*}

\address{Zi-Qiao Xu \\
Hebei Key Laboratory of Physics and Energy Technology\\
Department of Mathematics and Physics\\
North China Electric Power University \\
Yonghua Street 619, 071003, Baoding, P. R. China}
\email{ziqiaoxu678\symbol{64}gmail.com}

\address{Zhong-Xuan Mao \\
Hebei Key Laboratory of Physics and Energy Technology\\
Department of Mathematics and Physics\\
North China Electric Power University \\
Yonghua Street 619, 071003, Baoding, P. R. China}
\email{maozhongxuan2001\symbol{64}126.com}

\address{Jing-Feng Tian\\
Hebei Key Laboratory of Physics and Energy Technology\\
Department of Mathematics and Physics\\
North China Electric Power University\\
Yonghua Street 619, 071003, Baoding, P. R. China}
\email{tianjf\symbol{64}ncepu.edu.cn}

\subjclass[2000]{33C15, 33C05, 11B37, 41A05}
\keywords{Recurrence relation, Confluent hypergeometric function, Gauss hypergeometric function, Bessel function, Special function}
\thanks{*Corresponding author: Jing-Feng Tian, e-mail addresses,
tianjf@ncepu.edu.cn}

\begin{document}

\begin{abstract}
For  $a,b,c,z,p, \theta \in \mathbb{C}$, where $\mathbb{C}$ is the complex plane, $-c\notin \mathbb{N\cup }\left\{ 0\right\} $, let
\begin{equation*}
\mathcal{M}\left( z\right) =\left( 1-\theta z\right) ^{p}M\left(a;c;z\right) =\sum_{n=0}^{\infty }u_{n}z^{n},
\end{equation*}
where $|z| <\frac{1}{\theta}$,  $|\arg (1-\theta z)| < \pi$,
and let
\begin{equation*}
\mathcal{G}\left( z\right) =(1-\theta z) ^{p}F(a,b;c;z) =\sum_{n=0}^{\infty }v_{n} z^{n},
\end{equation*}
where $|z| < 1$,  $|\arg (1-\theta z)| < \pi$.
In this paper, we prove that the coefficients $u_{n}$ and $v_{n}$
for $n\geq 0$ satisfy a 3-order recurrence relation. These offer
a new way to study confluent hypergeometric function $M(a;c;z)$ and Gauss hypergeometric function $F(a,b;c;z)$. And we provide other special functions' recurrence relations of their coefficients, such as error function, Bessel function, incomplete gamma function, complete elliptic integral and Chebyshev polynomials.
\end{abstract}

\maketitle

\section{Introduction}
\numberwithin{equation}{section}
For complex numbers $a,c$, and $z$ with $c\neq 0,-1,-2,...$, the confluent hypergeometric function \cite{Buchholz-1969} is defined as:
\begin{equation*}
M(a;c;z) = \displaystyle\sum_{n=0}^{\infty}\frac{(a)_n}{(c)_n}\frac{z^n}{n!} = {}_1F_1(a;c;z),
\end{equation*}
where $\left( a\right) _{n}$ denotes
Pochhammer symbol defined by%
\begin{equation*}
\left( a\right) _{n}=a\left( a+1\right) \cdot \cdot \cdot \left(
a+n-1\right),
\end{equation*}%
for $n=1,2,...$, and $(a)_{0}=1$ for $a\neq 0$.
Confluent hypergeometric function is a solution of a confluent hypergeometric equation, which is a degenerate form of a hypergeometric differential equation:
\begin{equation*}
z\frac{d^2w}{dz^2} + (c-z)\frac{dw}{dz} - aw = 0
\end{equation*}
where two of the three regular singularities merge into an irregular singularity. As special case, the modified Bessel function and lower(upper) incomplete gamma functions $\gamma(s,z)$ ($\Gamma(s,z)$):
\begin{align*}
    & I_\alpha(x) = i^{-\alpha} J_\alpha(ix) = \sum_{m=0}^{\infty} \frac{1}{m! \Gamma(m + \alpha + 1)} \left(\frac{x}{2}\right)^{2m+\alpha}\\
    & \gamma(s,z) = \int_{0}^{z} t^{s-1}e^{-t}dt\\
    \intertext{and}
    & \Gamma(s,z) = \int_{z}^{\infty} t^{s-1}e^{-t}dt
\end{align*}
has the representations in confluent hypergeometric function:
\begin{equation*}
\begin{aligned}
    & I_{a-\frac{1}{2}}(z) = M(a;2a;z)\frac{1}{\Gamma(a+\frac{1}{2})}
    (\frac{z}{4})^{-\frac{1}{2} + a}e^\frac{-z}{2}
\end{aligned}
\end{equation*}
\begin{align*}
& \Gamma(s,z) = e^{-z}U(1-s;1-s;z)\\
\intertext{and}
& \gamma(s,z) = \frac{z^s}{s}M(1;s+1;-z),
\end{align*}
respectively, where $U(a;b;z) = \frac{\Gamma(1-b)}{\Gamma(a+1-b)}M(a;b;z) + \frac{\Gamma(b-1)}{\Gamma(a)}z^{1-b}M(a+1-b;2-b;z)$.

Also, Gauss hypergeometric function \cite{Beukers-2007,Aomoto-2011} is defined by
\begin{equation*}
F\left( a,b;c;z\right) =\sum_{n=0}^{\infty }\frac{\left( a\right) _{n}\left(
b\right) _{n}}{\left( c\right) _{n}}\frac{z^{n}}{n!}
\end{equation*}%
for complex numbers $a,b,c,z$ with $c\neq 0,-1,-2,...$ and $|z| < 1 $.
Gauss hypergeometric function is a solution of Euler's hypergeometric differential equation:
\begin{equation*}
z(1-z)\frac{d^2w}{dz^2} + [c-(a+b+1)z]\frac{dw}{dz} - abw = 0,
\end{equation*}
which has three regular singular points: $0,1$ and $\infty$, cf.\ \cite{Olver-2010}.
The complete elliptic integrals of the first and second kinds,
\[
K(r) = \int_0^{\pi/2} \frac{dt}{\sqrt{1 - r^2 \sin^2 t}}
\quad \text{and} \quad
E(r) = \int_0^{\pi/2} \sqrt{1 - r^2 \sin^2 t}\, dt,
\]
could be represented respectively by the hypergeometric function representations
\[
K(r) = \frac{\pi}{2} F\!\left( \frac{1}{2}, \frac{1}{2}; 1; r^2 \right)
\quad \text{and} \quad
E(r) = \frac{\pi}{2} F\!\left( -\frac{1}{2}, \frac{1}{2}; 1; r^2 \right).
\]

In the study of certain special functions, one often encounters the functions $\mathcal{M}(z)$ and $\mathcal{G}(z)$.
If $p$ is an integer, then it is easy to deal with $\mathcal{M}(z)$ and $\mathcal{G}(z)$. If $p$ is not an integer and $p \neq a + b - c$, then $\mathcal{M}(z)$ and $\mathcal{G}(z)$ are difficult to be treated. Then, how to deal with this kind of problem?

As is known to all, giving recurrence relations is one of important ways to study some special functions. In 2018, Yang and Tian \cite{Yang1} offered two recurrence relations of Maclaurin series coefficients of $(r')^pK(r)$ and $(r')^pE(r)$, $r,p \in \mathbb{R}$. In 2021, Yang and Tian \cite{YangandTian} presented another recurrence relation
of Maclaurin series coefficients of
$(1-x)^{-q}F(-\frac{1}{2},-\frac{1}{2};2;x)$, $x,q \in \mathbb{R}$, and provided the sharp lower and upper bounds for $E(r)$ in terms of weighted power means of arithmetic and geometric means. In 2025, Yang \cite{Yang2} provided the recurrence relation of the $(1-\theta x) ^{p}F(a,b;c;x)$, $\theta, p, a,b,c,x \in \mathbb{R}$.

Inspired by the recurrence formulas established in \cite{YangandTian} and \cite{Yang2}, the main aim of this paper is to establish a general recurrence relation of $\mathcal{M}(z)$ and $\mathcal{G}(z)$' Maclaurin series coefficients in complex fields.

\section{Recurrence relations of coefficients involving confluent hypergeometric function}

We start with a necessary lemma for the confluent hypergeometric funcion.
\begin{lemma}\label{lemma2.2}
For, $a,b,z \in \mathbb{C}$, $-b \notin \mathbb{N} \cup \{0\}$, let
\begin{equation*}
M \equiv M(a;c;z), \quad M(a-) \equiv M(a-1;c;z), \quad M(a+) \equiv M(a+1;c;z).
\end{equation*}
we have the following relation:
\begin{align*}
z \frac{dM}{dz} &= a(M(a+) - M) = (c-a)M(a-) + (a - c + z)M,
\end{align*}
while the first equation could also be written as:
\begin{align*}
    z \frac{dM(a-)}{dz} = (a-1)(M(a) - M(a-))
\end{align*}
\end{lemma}

\begin{theorem}\label{thm2.2}
Let $a, c, p, \theta, z \in \mathbb{C}$, $-c \notin \mathbb{N} \cup \{0\}$, and let $|\arg (1-\theta z)| < \pi$.
Suppose that
\[
    \mathcal{M}(z) = (1 - \theta z)^p M(a, c; z) = \sum_{n=0}^{\infty} u_n z^n, \quad |z| < \frac{1}{|\theta|}.
\]
Then $u_0 = 1$, $u_1 = \frac{a}{c} - \theta p$, $u_2 = \frac{1}{2} \left( \frac{a^2+a}{c^2+c} - \frac{2a\theta p}{c} + \theta^2(p-1)p \right)$ and
\[
    u_{n+1} = \beta_0(n)u_n + \beta_1(n)u_{n-1} + \beta_2(n)u_{n-2},
\]
where
\begin{align*}
\beta_0(n) &= \frac{a + 2\theta n(c+n-1) - \theta p(c+2n) + n}{(n+1)(c+n)}, \\
\beta_1(n) &= \frac{\theta(-2a - \theta(n-p-1)(c+n-p-2) - 2n + p + 2)}{(n+1)(c+n)}, \\
\beta_2(n) &= \frac{\theta^2(a+n-p-2)}{(n+1)(c+n)}.
\end{align*}

\begin{proof}
Let $\mathcal{M}^{*}(z) = (1 - \theta z)^{p}M(a-1;c;z) =
(1 - \theta z)^{p}M(a-) = \sum_{n=0}^{\infty} u'_n z^n$.
We first differentiate both two power series expressions and multiply $z(1 - \theta z)$, then we get
\begin{equation*}
    -p \theta z\mathcal{M}(z) + (1- \theta z)((c-a) \mathcal{M}^{*}(z) + (a-c+z)\mathcal{M}(z)) =
    (z - \theta z^2)\displaystyle\sum_{n=0}^{\infty}
    nu_nz^{n-1},
\end{equation*}
and
\begin{equation*}
    -p \theta z\mathcal{M}^{*}(z) + (a-1)(1- \theta z)(\mathcal{M}(z) - \mathcal{M}^{*}(z)) =
    (z - \theta z^2)\displaystyle\sum_{n=0}^{\infty}
    nu_n'z^{n-1}.
\end{equation*}

Then, we just expand $\mathcal{M}$ and $\mathcal{M}^{*}$ into power series:
\begin{equation*}
\begin{aligned}
    &(b-a)u'_0 + (a-c)u_0 + z[-p \theta u_0 + (c-a)u_1' - \theta (b-c)u_0' + (a-c)u_1 + (1-\theta (a-c))u_0]
    + \\
    &\displaystyle\sum_{n=2}^{\infty} z^n(-p \theta u_{n-1} + (c-a)u_n' - \theta (c-a)u_{n-1}' + (a-c)u_n + (1-\theta (a-c))u_{n-1} - \theta u_{n-2})\\
    &= u_1z + \displaystyle\sum_{n=0}^{\infty}(nu_n -
    \theta u_{n-1})z^n
\end{aligned}
\end{equation*}
and
\begin{equation*}
\begin{aligned}
    &(a-1)(u_0 - u_0') + z(-p \theta u_0' + (a-1)u_1 - (a-1)\theta u_0 - (a-1)u_1' + (a-1)\theta u_0')
    + \\
    &\displaystyle\sum_{n=2}^{\infty}z^n(-p \theta u_{n-1}' + (a-1)(u_n - u_n') - \theta (a-1)(u_{n-1} - u_{n-1}'))\\
    & = zu_1' + \displaystyle\sum_{n=2}^{\infty} z^n(nu_n' - (n-1)u_{n-1}').
\end{aligned}
\end{equation*}

Compare the coefficients of $z^n$ between these equations, then we get
\begin{equation*}
\begin{aligned}
    & u_0 = 1\\
    & u_1 = \frac{a}{c} - \theta p\\
    & (c-a)(u_n' - u_{n-1}') = (1 + c - a)u_n +
    (\theta (a-c+1) + p\theta - 1)u_{n-1} + \theta u_{n-2}\\
    & (a - 1)(u_n - u_{n-1}) = au_n' + (p\theta -\theta a)u_{n-1}'
\end{aligned}
\end{equation*}
Represent the $u_{n-1}'$ by $u_n$, $u_{n-1}$ and $u_{n-2}$, we obtain
\begin{equation*}
u'_{n-1} = \frac{C_0 u_n + C_1 u_{n-1} + C_2 u_{n-2}}{\theta (a - c) (p + 1)},
\end{equation*}
where, $\begin{alignedat}[t]{2}
C_0 &= n(c+n-1), \qquad
& C_1 &= \theta\bigl(ap+a-cn-n^2+np+2n-p-1\bigr)-a-n+1, \\
C_2 &= \theta(a+n-1).
\end{alignedat}$

Now we substitute the recurrence relation of $u_n'$ into the original equation, and finally get the recurrence relation of $u_n$. While for $u_2$, we just consider the Taylor expansion of $(1-\theta z)^p$ and $M(a;c;z)$:
\begin{equation*}
\begin{aligned}
    & (1-\theta z)^p = 1 - p\theta + \frac{p(p-1)\theta }{2}z^2 + o(z^3)\\
    & M(a;c;z) = 1 + \frac{a}{c} z + \frac{a(a+1)}{c(c+1)}z^2 + o(z^3)
\end{aligned}
\end{equation*}
It is easy to see that the coefficient of $z^2$ is
$u_2 = \frac{1}{2} \left( \frac{a^2+a}{c^2+c} - \frac{2a\theta p}{c} + \theta^2(p-1)p \right)$,
which completes the proof.
\end{proof}
\end{theorem}

If we restrict $p \in \mathbb{Z}$, we would get the following result from Theorem \ref{thm2.2} without setting $\arg (1-\theta z) \in (-\pi,\pi)$ canonically.
\begin{theorem}\label{p in Z}
Let $a, c, \theta, z \in \mathbb{C}$, $p \in \mathbb{Z}$, $-c \notin \mathbb{N} \cup \{0\}$. Then
\[
    (1 - \theta z)^p M(a, c; z) = \sum_{n=0}^{\infty} u_n z^n, \quad |z| < \frac{1}{|\theta|},
\]
with $u_0 = 1$, $u_1 = \frac{a}{c} - \theta p$, $u_2 = \frac{1}{2} \left( \frac{a^2+a}{c^2+c} - \frac{2a\theta p}{c} + \theta^2(p-1)p \right)$ and
\[
    u_{n+1} = \beta_0(n)u_n + \beta_1(n)u_{n-1} + \beta_2(n)u_{n-2},
\]
where
\begin{align*}
\beta_0(n) &= \frac{a + 2\theta n(c+n-1) - \theta p(c+2n) + n}{(n+1)(c+n)}, \\
\beta_1(n) &= \frac{\theta(-2a - \theta(n-p-1)(c+n-p-2) - 2n + p + 2)}{(n+1)(c+n)}, \\
\beta_2(n) &= \frac{\theta^2(a+n-p-2)}{(n+1)(c+n)}.
\end{align*}
\end{theorem}

%---corollary for theta--------
Now we will give some special case of Theorem \ref{thm2.2} for $\theta = 1,-1,i,1+i$.
Substitute $\theta = 1$ into Theorem \ref{thm2.2}, we could get the following result.
\begin{corollary}
Let $a, c, z, p \in \mathbb{C}$, $-c \notin \mathbb{N} \cup \{0\}$, and let $|\arg (1-z)| < \pi$. Then
\[
    (1 -  z)^p M(a, c; z) = \sum_{n=0}^{\infty} u_n z^n, \quad |z| < 1,
\]
with $u_0 = 1$, $u_1 = \frac{a}{c} -  p$, $u_2 = \frac{1}{2} \left( \frac{a^2+a}{c^2+c} - \frac{2ap}{c} + (p-1)p \right)$ and
\[
    u_{n+1} = \beta_0(n)u_n + \beta_1(n)u_{n-1} + \beta_2(n)u_{n-2},
\]
where
\begin{align*}
\beta_0(n) &= \frac{a + 2 n(c+n-1) -  p(c+2n) + n}{(n+1)(c+n)}, \\
\beta_1(n) &= \frac{(-2a - (n-p-1)(c+n-p-2) - 2n + p + 2)}{(n+1)(c+n)}, \\
\beta_2(n) &= \frac{(a+n-p-2)}{(n+1)(c+n)}.
\end{align*}
\end{corollary}

Substitute $\theta = -1$ into Theorem \ref{thm2.2}, we could get the following result.
\begin{corollary}
Let $a, c,z,p  \in \mathbb{C}$, $-c \notin \mathbb{N} \cup \{0\}$, and let $|\arg (1+z)| < \pi$. Then
\[
    (1 + z)^p M(a, c; z) = \sum_{n=0}^{\infty} u_n z^n, \quad |z| < 1,
\]
with $u_0 = 1$, $u_1 = \frac{a}{c} + p$, $u_2 = \frac{1}{2} \left( \frac{a^2+a}{c^2+c} + \frac{2a p}{c} + (p-1)p \right)$ and
\[
    u_{n+1} = \beta_0(n)u_n + \beta_1(n)u_{n-1} + \beta_2(n)u_{n-2},
\]
where
\begin{align*}
\beta_0(n) &= \frac{a - 2 n(c+n-1) + p(c+2n) + n}{(n+1)(c+n)}, \\
\beta_1(n) &= \frac{(2a - (n-p-1)(c+n-p-2) + 2n - p - 2)}{(n+1)(c+n)}, \\
\beta_2(n) &= \frac{(a+n-p-2)}{(n+1)(c+n)}.
\end{align*}
\end{corollary}

Substitute $\theta = i$ into Theorem \ref{thm2.2}, then we have the following result.
\begin{corollary}
Let $a, c,z,p \in \mathbb{C}$, $-c \notin \mathbb{N} \cup \{0\}$, and let $|\arg (1-iz)| < \pi$. Then
\[
    (1 - i z)^p M(a, c; z) = \sum_{n=0}^{\infty} u_n z^n, \quad |z| < 1,
\]
with $u_0 = 1$, $u_1 = \frac{a}{c} - i p$, $u_2 = \frac{1}{2} \left( \frac{a^2+a}{c^2+c} - \frac{2ia p}{c} - (p-1)p \right)$ and
\[
    u_{n+1} = \beta_0(n)u_n + \beta_1(n)u_{n-1} + \beta_2(n)u_{n-2},
\]
where
\begin{align*}
\beta_0(n) &= \frac{a + 2i n(c+n-1) - i p(c+2n) + n}{(n+1)(c+n)}, \\
\beta_1(n) &= \frac{i(-2a - i(n-p-1)(c+n-p-2) - 2n + p + 2)}{(n+1)(c+n)}, \\
\beta_2(n) &= -\frac{(a+n-p-2)}{(n+1)(c+n)}.
\end{align*}
\end{corollary}

Substitute $\theta = 1+i$ into Theorem \ref{thm2.2}, we could get the following result.
\begin{corollary}
Let $a, c,z,p \in \mathbb{C}$, $-c \notin \mathbb{N} \cup \{0\}$, and let $|\arg (1-(1+i)z)| < \pi$. Then
\[
    (1 - (1+i) z)^p M(a, c; z) = \sum_{n=0}^{\infty} u_n z^n, \quad |z| < \frac{1}{\sqrt{2}},
\]
with $u_0 = 1$, $u_1 = \frac{a}{c} - (1+i) p$, $u_2 = \frac{1}{2} \left( \frac{a^2+a}{c^2+c} - \frac{2a(1+i) p}{c} + 2i(p-1)p \right)$ and
\[
    u_{n+1} = \beta_0(n)u_n + \beta_1(n)u_{n-1} + \beta_2(n)u_{n-2},
\]
where
\begin{align*}
\beta_0(n) &= \frac{a + 2(1+i) n(c+n-1) - (1+i) p(c+2n) + n}{(n+1)(c+n)}, \\
\beta_1(n) &= \frac{(1+i)(-2a - (1+i)(n-p-1)(c+n-p-2) - 2n + p + 2)}{(n+1)(c+n)}, \\
\beta_2(n) &= \frac{2i(a+n-p-2)}{(n+1)(c+n)}.
\end{align*}
\end{corollary}

\section{Recurrence relations of coefficients involving gauss hypergeometric function}
\numberwithin{equation}{section}
We start with a useful lemma for Gauss hypergeometric function.

\begin{lemma}\label{L-dF,F+}
\cite{Yang2} Let $a,b,c,z \in \mathbb{C}$,$-c\notin \mathbb{N\cup }\left\{ 0\right\} $ and $|z|<1$
\begin{equation*}
\begin{array}{ccccc}
F\equiv F\left( a,b;c,z\right) , &  & F_{a-}\equiv F\left( a-1,b;c,z\right),
&  & F_{a+}\equiv F\left( a+1,b;c,z\right)%.
\end{array}%
\end{equation*}
Then
\begin{eqnarray}
\frac{dF}{dz} &=&\frac{\left( c-a\right) F_{a-}+\left( a-c+bz\right) F}{%
z\left( 1-z\right) },  \label{dF} \\
\frac{dF_{a-}}{dz} &=&\frac{a-1}{z}\left( F-F_{a-}\right)  \label{dF-}
\end{eqnarray}
\end{lemma}

We are now going to present and prove the main theorem for Gauss hypergeometric function.

\begin{theorem}\label{main thm1}
\label{T-Ut}Let $a,b,c,p, \theta,z \in \mathbb{C}$, $-c\notin \mathbb{N\cup }\left\{ 0\right\} $. Then we have
\begin{equation*}
\mathcal{G}(z) =\left( 1-\theta z\right) ^{p}F\left(
a,b;c;z\right) =\sum_{n=0}^{\infty }v_{n}z^{n},
\end{equation*}
with $|z| < 1$, $|\arg (1-\theta z)| < \pi$, $v_{0}=1$, $v_{1}=ab/c-p\theta $,%
\begin{equation*}
v_{2}=\frac{1}{2}\theta ^{2}p\left( p-1\right) -\theta p\frac{ab}{c}+\frac{1%
}{2}\frac{ab\left( b+1\right) \left( a+1\right) }{c\left( c+1\right) },
\end{equation*}
and for $n\geq 2$,
\begin{equation}
v_{n+1}=\frac{\mu _{n,p,\theta }\left( a,b,c\right) }{\left( n+1\right)
\left( n+c\right) }v_{n}-\theta \dfrac{\zeta _{n,p,\theta }\left(
a,b,c\right) }{\left( n+1\right) \left( n+c\right) }v_{n-1}+\theta ^{2}\frac{%
\sigma _{n,p}\left( a,b\right) }{\left( n+1\right) \left( n+c\right) }%
v_{n-2},  \label{un+1-n-2-rr}
\end{equation}%
where
\begin{eqnarray*}
\mu _{n,p,\theta }\left( a,b,c\right) &=&\left( n+a\right) \left( n+b\right)
+\theta \left[ 2n^{2}-2n\left( p-c+1\right) -cp\right] , \\
\zeta _{n,p,\theta }\left( a,b,c\right) &=&2n^{2}+2\left( a+b-p-2\right)
n-\left( a+b-1\right) p+2\left( a-1\right) \left( b-1\right) \\
&&+\theta \left( n-p-1\right) \left( n-p+c-2\right) , \\
\sigma _{n,p}\left( a,b\right) &=&\left( n+a-p-2\right) \left(
n+b-p-2\right).
\end{eqnarray*}
\end{theorem}

\begin{proof}
Denote
\begin{equation*}
\mathcal{G}^{\ast}(z) \equiv \left( 1-\theta z\right) ^{p}F\left(
a-1,b;c;z\right) =\sum_{n=0}^{\infty }v_{n}^{\ast }z^{n}.  \label{Ut*}
\end{equation*}

The proof is analogous to that of Theorem \ref{thm2.2} in Section 2, we first consider differentiating $\mathcal{G}(z)$
\begin{equation*}
-p\theta \left( 1-\theta z\right) ^{p-1}F+\left( 1-\theta z\right)
^{p}F^{\prime }=\sum_{n=0}^{\infty }nv_{n}z^{n-1}.
\end{equation*}%
Substituting (\ref{dF}) into the above equation and multiplying by $z\left( 1-z\right) \left( 1-\theta z\right) $ give
\begin{equation}
\begin{array}{c}
-p\theta z\left( 1-z\right) \mathcal{G}(z) +\left( c-a\right)
\left( 1-\theta z\right) \mathcal{G}^{\ast}(z) +\left(
a-c+bz\right) \left( 1-\theta z\right) \mathcal{G}(z)
\\
=\left( 1-z\right) \left( 1-\theta z\right) \sum_{n=0}^{\infty
}nv_{n}z^{n}.
\end{array}
\label{Ut-Ut*}
\end{equation}%
Now we consider expanding $\mathcal{G}$ and $\mathcal{G}^*$ into the assumed power series:
\begin{eqnarray*}
&&-p\theta \left( \sum_{n=1}^{\infty }v_{n-1}z^{n}-\sum_{n=2}^{\infty
}v_{n-2}z^{n}\right) +\left( c-a\right) \left( \sum_{n=0}^{\infty
}v_{n}^{\ast }z^{n}-\theta \sum_{n=1}^{\infty }v_{n-1}^{\ast }z^{n}\right) \\
&&+\left( a-c\right) \sum_{n=0}^{\infty }v_{n}z^{n}+\left( b-\theta \left(
a-c\right) \right) \sum_{n=1}^{\infty }v_{n-1}z^{n}-b\theta
\sum_{n=2}^{\infty }v_{n-2}z^{n} \\
&=&\sum_{n=1}^{\infty }nv_{n}z^{n}-\left( \theta +1\right)
\sum_{n=2}^{\infty }\left( n-1\right) v_{n-1}z^{n}+\theta \sum_{n=2}^{\infty
}\left( n-2\right) v_{n-2}z^{n}.
\end{eqnarray*}%
We may rewrite the above equation like
\begin{eqnarray*}
&&\left( c-a\right) \left( v_{0}^{\ast }-v_{0}\right) +\left[ -p\theta
v_{0}+\left( c-a\right) \left( v_{1}^{\ast }-\theta v_{0}^{\ast }\right)
+\left( a-c\right) v_{1}+\left( b-\theta \left( a-c\right) \right) v_{0}%
\right] z \\
&&+\sum_{n=2}^{\infty }\left[ \left( c-a\right) \left( v_{n}^{\ast }-\theta
v_{n-1}^{\ast }\right) +\left( a-c\right) v_{n}+\left( b-\theta \left(
a-c+p\right) \right) v_{n-1}+\left( p-b\right) \theta v_{n-2}\right] z^{n} \\
&=&v_{1}z+\sum_{n=2}^{\infty }\left[ nv_{n}-\left( \theta +1\right) \left(
n-1\right) v_{n-1}+\theta \left( n-2\right) v_{n-2}\right] z^{n}.
\end{eqnarray*}
Then, we could give the following results by comparing coefficients of $z^{n}$,
\begin{equation*}
\left( c-a\right) \left(
v_{0}^{\ast }-v_{0}\right) =0, \quad
\left( c-a\right) \left( v_{1}^{\ast }-\theta v_{0}^{\ast }\right) +\left(
b-\theta \left( a-c+p\right) \right) v_{0}+\left( a-c\right) v_{1}=v_{1},
\end{equation*}
\begin{equation}\label{un*-un-r1}
\begin{aligned}
&\left( c-a\right) \left( v_{n}^{\ast }-\theta v_{n-1}^{\ast }\right) = \left(
n-a+c\right) v_{n}+\theta \left( n-2+b-p\right) v_{n-2} \\
&+\left[ \theta \left( a-c+p\right) -b-\left( \theta +1\right) \left(
n-1\right) \right] v_{n-1}\text{ for }n\geq 2.
\end{aligned}
\end{equation}

In the same way, we will differentiate $\mathcal{G}^{\ast }$ and get
\begin{equation*}
-p\theta \left( 1-\theta z\right) ^{p-1}F_{a-}+\left( 1-\theta z\right)
^{p}F_{a-}^{\prime }=\sum_{n=0}^{\infty }nv_{n}^{\ast }z^{n-1}.
\end{equation*}%

Substituting (\ref{dF-}) into the above equation and multiplying by $x\left( 1-\theta x\right) $ give:
\begin{equation*}
-p\theta z\mathcal{G}^{\ast}(z) +\left( a-1\right) \left(
1-\theta z\right) \left( \mathcal{G}(z) -\mathcal{G}^{\ast
}(z) \right) =\left( 1-\theta z\right) \sum_{n=0}^{\infty
}nv_{n}^{\ast }z^{n}.
\end{equation*}%
Then, we still consider expanding $\mathcal{G}$ and $\mathcal{G}^*$ in power series and get
\begin{equation*}
-p\theta \sum_{n=1}^{\infty }v_{n-1}^{\ast }z^{n}+\left( a-1\right) \left(
1-\theta z\right) \sum_{n=0}^{\infty }\left( v_{n}-v_{n}^{\ast }\right)
z^{n}=\left( 1-\theta z\right) \sum_{n=0}^{\infty }nv_{n}^{\ast }z^{n},
\end{equation*}%
which can be rewritten as
\begin{eqnarray*}
&&\left( a-1\right) \left( v_{0}-v_{0}^{\ast }\right) -p\theta
\sum_{n=1}^{\infty }v_{n-1}^{\ast }z^{n}+\left( a-1\right)
\sum_{n=1}^{\infty }\left( v_{n}-\theta v_{n-1}-v_{n}^{\ast }+\theta
v_{n-1}^{\ast }\right) z^{n} \\
&=&\sum_{n=1}^{\infty }\left( nv_{n}^{\ast }-\theta \left( n-1\right)
v_{n-1}^{\ast }\right) z^{n}.
\end{eqnarray*}%
From the equality of the coefficients of $z^{n}$, it follows that $\left( a-1\right) \left(
v_{0}-v_{0}^{\ast }\right) =0$ and%
\begin{equation}\label{un*-un-r2}
\left( n+a-1\right) v_{n}^{\ast }-\theta \left( n-2+a-p\right) v_{n-1}^{\ast
}=\left( a-1\right) v_{n}-\theta \left( a-1\right) v_{n-1}\text{ for }n\geq 1%
\text{.}
\end{equation}

Now from equations (\ref{un*-un-r1}) and (\ref{un*-un-r2}), we obtain the representation for $v_{n-1}^{\ast}$ by $v_n,v_{n-1},v_{n-2}$ with $n\geq 2$
\begin{equation}
\begin{array}{c}
\theta \left( p+1\right) \left( a-c\right) v_{n-1}^{\ast }=n\left(
n+c-1\right) v_{n}+\theta \left( n+a-1\right) \left( n+b-p-2\right)
v_{n-2}\\
-\left[ \left( \theta +1\right) n^{2}+\left( a+b+\left( c-p-2\right) \theta
-2\right) n+\left( a-1\right) \left( b-\left( p+1\right) \theta -1\right)
\right] v_{n-1}.
\end{array}
\label{un*-un-n-2}
\end{equation}%
When $\theta \left( p+1\right) \neq 0$, substituting the above expressions of $v_{n-1}^{\ast }$ and $v_{n}^{\ast }$
into (\ref{un*-un-r1}) gives
\begin{equation*}
\begin{array}{l}
v_{n+1}-\dfrac{\left( 2\theta +1\right) n^{2}+\left( a+b+2\theta \left(
c-1-p\right) \right) n+\left( ab-cp\theta \right) }{\left( n+1\right) \left(
n+c\right) }v_{n} \\
+\theta \dfrac{
\begin{array}{c}
\left( \theta +2\right) n^{2}+\left( 2a+2b-2p+\theta \left( c-2p-3\right)
-4\right) n \\
+\theta \left( p+1\right) \left( p-c+2\right) +2ab-\left( p+2\right) \left(
a+b-1\right)
\end{array}
}{\left( n+1\right) \left( n+c\right) }v_{n-1} \\
\end{array}
\label{un+1-n-2-rra}
\end{equation*}
\begin{equation*}
    -\theta ^{2}\dfrac{\left( n+a-p-2\right) \left( n+b-p-2\right) }{\left(
n+1\right) \left( n+c\right) }v_{n-2}=0\text{ for }n\geq 2,
\end{equation*}

which implies (\ref{un+1-n-2-rr}).

Also, we have the following expression by using the Cauchy product formula
\begin{equation*}
v_{n}=\sum_{k=0}^{n}\frac{\left( a\right) _{k}\left( b\right) _{k}}{k!\left(
c\right) _{k}}\frac{\theta ^{n-k}\left( -p\right) _{n-k}}{\left( n-k\right) !%
},
\end{equation*}%
which means $v_{0}=1$,%
\begin{equation*}
v_{1}=\frac{ab}{c}-p\theta \text{ \ and \ }v_{2}=\frac{1}{2}\theta
^{2}p\left( p-1\right) -\theta p\frac{ab}{c}+\frac{1}{2}\frac{ab\left(
b+1\right) \left( a+1\right) }{c\left( c+1\right) }.
\end{equation*}

When $\theta =0$, the recurrence relation (\ref{un+1-n-2-rr}) holds trivially. When $p=-1$, it could be seen that
\begin{equation*}
v_{n}=\sum_{k=0}^{n}\frac{\left( a\right) _{k}\left( b\right) _{k}}{k!\left(
c\right) _{k}}\theta ^{n-k},
\end{equation*}%
which could be written as
\begin{equation*}
v_{n+1}-\theta v_{n}=\frac{\left( a\right) _{n+1}\left( b\right) _{n+1}}{%
\left( n+1\right) !\left( c\right) _{n+1}}.
\end{equation*}%
This also satisfied the recurrence relation (\ref{un+1-n-2-rr}), and the whole proof is done.
\end{proof}

%----collary-----------------
We have the following result by substituting $\theta = -1$ into Theorem \ref{main thm1}.
\begin{corollary}\label{col2.1}
\label{C-th=-1}Let $a,b,c,p,z \in \mathbb{C}$, $-c\notin \mathbb{N\cup }\left\{ 0\right\} $, $|z| < 1$, $|\arg (1+z)| < \pi$. Then
\begin{equation*}
\left( 1+z\right) ^{p}F\left( a,b;c;z\right)
=\sum_{n=0}^{\infty }v_{n}z^{n},
\end{equation*}
with $v_{0}=1$, $v_{1}=ab/c+p$,
\begin{equation*}
v_{2}=\frac{1}{2}p\left( p-1\right) +p\frac{ab}{c}+\frac{1}{2}\frac{ab\left(
b+1\right) \left( a+1\right) }{c\left( c+1\right) },
\end{equation*}%
and for $n\geq 2$,%
\begin{equation*}
v_{n+1}=\frac{\mu _{n,p}\left( a,b,c\right) }{\left( n+1\right) \left(
n+c\right) }v_{n}+\dfrac{\zeta _{n,p}\left( a,b,c\right) }{\left( n+1\right)
\left( n+c\right) }v_{n-1}+\frac{\sigma _{n,p}\left( a,b\right) }{\left(
n+1\right) \left( n+c\right) }v_{n-2},
\end{equation*}%
where
\begin{eqnarray*}
\mu _{n,p}\left( a,b,c\right) &=&-n^{2}+\left( a+b-2c+2p+2\right) n+\left(
ab+cp\right) , \\
\zeta _{n,p}\left( a,b,c\right) &=&\left( n+2a+2b-c\right) \left( n-1\right)
-p^{2}-\left( a+b-c+2\right) p+2ab, \\
\sigma _{n,p}\left( a,b\right) &=&\left( n+a-p-2\right) \left(
n+b-p-2\right) .
\end{eqnarray*}
\end{corollary}

Substituting $\theta =1$ into Theorem \ref{main thm1}, we can easily obtain a special 3-order recurrence relation. However, the order can be reduced from 3 to 2 from the proof of Theorem \ref{T-Ut}.

\begin{corollary}\label{col2.3}
\label{C-th=1}Let For $a,b,c,p,z\in \mathbb{C}$, $-c\notin \mathbb{N\cup }\left\{ 0\right\} $. Then
\begin{equation*}
\mathcal{G}(z) = \left( 1-z\right) ^{p}F\left( a,b;c;z\right)
=\sum_{n=0}^{\infty }u_{n}z^{n},
\end{equation*}%
with $|z| < 1$, $|\arg (1-z)| < \pi$, $v_{0}=1$, $v_{1}=ab/c-p$ and for $n\geq 1$,%
\begin{equation}
v_{n+1}=\tau_{n}v_{n}-\rho_{n}v_{n-1},  \label{u-th=1-rr}
\end{equation}%
where
\begin{eqnarray}
\tau _{n} &=&\dfrac{2n^{2}+\left( a+b+c-2p-1\right) n+ab-cp}{%
\left( n+1\right) \left( n+c\right) },  \label{aln} \\
\rho _{n} &=&\dfrac{\left( n+a-p-1\right) \left( n+b-p-1\right) }{\left(
n+1\right) \left( n+c\right) }.  \label{ben}
\end{eqnarray}
\end{corollary}

\begin{proof}
When $\theta =1$, dividing $\left( 1-z\right) $ by the relation (\ref{Ut-Ut*}%
) gives%
\begin{equation*}
-pz\mathcal{G}(z) +\left( c-a\right) \mathcal{G}^{\ast }(z)
+\left( a-c+bz\right) \mathcal{G}(z) =\left( 1-z\right)
\sum_{n=0}^{\infty }nv_{n}z^{n}.
\end{equation*}%
Expanding $\mathcal{G}(z)$ and $\mathcal{G}^{\ast }(z)$ in power series gives
\begin{eqnarray*}
&&-p\sum_{n=1}^{\infty }v_{n-1}z^{n}+\left( c-a\right) \sum_{n=0}^{\infty
}v_{n}^{\ast }z^{n}+\left( a-c\right) \sum_{n=0}^{\infty
}v_{n}z^{n}+b\sum_{n=1}^{\infty }v_{n-1}z^{n} \\
&=&\sum_{n=1}^{\infty }nv_{n}z^{n}-\sum_{n=1}^{\infty }\left( n-1\right)
v_{n-1}z^{n},
\end{eqnarray*}%
Comparing coefficients of $z^{n}$ gives $\left( c-a\right) \left(
v_{0}^{\ast }-v_{0}\right) =0$,%
\begin{equation*}
-pv_{n-1}+\left( c-a\right) v_{n}^{\ast }+\left( a-c\right)
v_{n}+bv_{n-1}=nv_{n}-\left( n-1\right) v_{n-1}\text{ for }n\geq 1\text{,}
\end{equation*}%
which implies%
\begin{equation}
\left( c-a\right) v_{n}^{\ast }=\left( n-a+c\right) v_{n}-\left(
n-1+b-p\right) v_{n-1}\text{ for }n\geq 1\text{.}  \label{un*-un,n-1}
\end{equation}
When $c\neq a$, eliminating $v_{n}^{\ast }$ and $v_{n-1}^{\ast }$ from the
Eqs. (\ref{un*-un,n-1}) and (\ref{un*-un-r2}) for $\theta =1$ yields%
\begin{eqnarray*}
v_{n} &=&\frac{2\left( n-1\right) ^{2}+\left( a+b+c-2p-1\right) \left(
n-1\right) +ab-cp}{n\left( c+n-1\right) }v_{n-1} \\
&&-\frac{\left( n+a-p-2\right) \left( n+b-p-2\right) }{n\left( c+n-1\right) }%
v_{n-2}\text{ for }n\geq 2.
\end{eqnarray*}%
Replacing $n$ by $n+1$ give the desired recurrence formula (\ref{u-th=1-rr}).

When $c=a$, from the relation (\ref{un*-un,n-1}) we have%
\begin{equation*}
v_{n}=\frac{n-1+b-p}{n}v_{n-1}\text{ for }n\geq 1,
\end{equation*}%
which, by an easy check, satisfies the recurrence relation (\ref{u-th=1-rr}).

The values of $v_{0}$ and $v_{1}$ follows Theorem \ref{T-Ut}, and the proof
is done.
\end{proof}

We have the following result by substituting $\theta = -i$ into Theorem \ref{main thm1}.
\begin{corollary}\label{col2.4}
Let $a,b,c,p,z \in \mathbb{C}$, $-c\notin \mathbb{N\cup }\left\{ 0\right\}$. Then
\begin{equation*}
\left( 1-iz\right) ^{p}F\left( a,b;c;z\right)
=\sum_{n=0}^{\infty }u_{n}z^{n},
\end{equation*}%
with $|z| < 1$, $|\arg (1-iz)| < \pi$, $v_{0}=1$, $v_{1}=ab/c - ip$,%
\begin{equation*}
v_{2}=-\frac{1}{2}p\left( p-1\right) - ip\frac{ab}{c}+\frac{1}{2}\frac{ab\left(
b+1\right) \left( a+1\right) }{c\left( c+1\right) },
\end{equation*}%
and for $n\geq 2$,%
\begin{equation*}
v_{n+1}=\frac{\mu _{n,p}\left( a,b,c\right) }{\left( n+1\right) \left(
n+c\right) }v_{n} - i\dfrac{\zeta _{n,p}\left( a,b,c\right) }{\left( n+1\right)
\left( n+c\right) }v_{n-1} - \frac{\sigma _{n,p}\left( a,b\right) }{\left(
n+1\right) \left( n+c\right) }v_{n-2},
\end{equation*}%
where
\begin{eqnarray*}
\mu _{n,p}\left( a,b,c\right) &=&\left( n+a\right) \left( n+b\right)
+i \left[ 2n^{2}-2n\left( p-c+1\right) -cp\right] , \\
\zeta _{n,p}\left( a,b,c\right) &=&2n^{2}+2\left( a+b-p-2\right)
n-\left( a+b-1\right) p+2\left( a-1\right) \left( b-1\right) \\
&&+i \left( n-p-1\right) \left( n-p+c-2\right) , \\
\sigma _{n,p}\left( a,b\right) &=&\left( n+a-p-2\right) \left(
n+b-p-2\right) .
\end{eqnarray*}
\end{corollary}

 We have the following result by substituting $\theta = 1+i$ into Theorem \ref{main thm1}.
\begin{corollary}\label{col2.5}
Let $a,b,c,p,z \in \mathbb{C}$, $-c\notin \mathbb{N\cup }\left\{ 0\right\} $. Then
\begin{equation*}
\left( 1-(1+i)z\right) ^{p}F\left( a,b;c;z\right)
=\sum_{n=0}^{\infty }u_{n}z^{n},
\end{equation*}%
with $|z| < 1$, $z \neq \frac{1}{1+i}$, $|\arg (1-(1+i)z)| < \pi$, $v_{0}=1$, $v_{1}=ab/c - (1+i)p$,%
\begin{equation*}
v_{2}=-\frac{1}{2}p\left( p-1\right) - (1+i)p\frac{ab}{c}+\frac{1}{2}\frac{ab\left(
b+1\right) \left( a+1\right) }{c\left( c+1\right) },
\end{equation*}%
and for $n\geq 2$,%
\begin{equation*}
v_{n+1}=\frac{\mu _{n,p}\left( a,b,c\right) }{\left( n+1\right) \left(
n+c\right) }v_{n} - (1+i)\dfrac{\zeta _{n,p}\left( a,b,c\right) }{\left( n+1\right)
\left( n+c\right) }v_{n-1} + 2i \frac{\sigma _{n,p}\left( a,b\right) }{\left(
n+1\right) \left( n+c\right) }v_{n-2},
\end{equation*}%
where
\begin{eqnarray*}
\mu _{n,p}\left( a,b,c\right) &=&\left( n+a\right) \left( n+b\right)
+(1+i) \left[ 2n^{2}-2n\left( p-c+1\right) -cp\right] , \\
\zeta _{n,p}\left( a,b,c\right) &=&2n^{2}+2\left( a+b-p-2\right)
n-\left( a+b-1\right) p+2\left( a-1\right) \left( b-1\right) \\
&&+(1+i) \left( n-p-1\right) \left( n-p+c-2\right) , \\
\sigma _{n,p}\left( a,b\right) &=&\left( n+a-p-2\right) \left(
n+b-p-2\right) .
\end{eqnarray*}
\end{corollary}

\section{Applications}
We are now going to use the main results in Section 2 and 3 to give some Maclaurin recurrence relations involving other special functions.

%------corollary for confluent--------

\subsection{Applications of the confluent hypergeometric function}
\mbox{} % 这里加一个空盒子并强制换行

First, we present a recurrence relation involving the error function $Erf(z)$.

\begin{corollary}
For $z,p \in \mathbb{C}$, $\quad |z|^2 < \frac{1}{|\theta|}$ and $\arg(1+\theta z^2) \in (-\pi,\pi)$, we could have a recurrence relation of coefficients for the error function:
\begin{equation*}
(1+\theta z^2)^p Erf(z) = \sum_{n=0}^{\infty}a_nz^{2n+1}
\end{equation*}
with $a_n = \frac{(-1)^nu_n}{\sqrt{\pi}}$, $u_0 = 1$, $u_1 = \frac{1}{3} - \theta p$, $u_2 = \frac{1}{2} \left( \frac{1}{5} - \frac{2\theta p}{3} + \theta^2(p-1)p \right)$ and
\[
    u_{n+1} = \beta_0(n)u_n + \beta_1(n)u_{n-1} + \beta_2(n)u_{n-2},
\]
where
\begin{align*}
\beta_0(n) &= \frac{1 + 2\theta n(1+2n) - \theta p(3+4n) + n}{(n+1)(3+2n)}, \\
\beta_1(n) &= \frac{\theta(-1 - \theta(n-p-1)(\frac{3}{2}+n-p-2) - 2n + p + 2)}{(n+1)(\frac{3}{2}+n)}, \\
\beta_2(n) &= \frac{\theta^2(2n-2p-3)}{(n+1)(3+2n)},
\end{align*}

\begin{proof}
It is known all that:
\begin{equation*}
Erf(z) = \frac{z}{\sqrt{\pi}}M(\frac{1}{2};\frac{3}{2};-z^2)
\end{equation*}
Then, applying Theorem \ref{thm2.2}, we have:
\begin{equation*}
(1+\theta z^2)^p Erf(z) = \frac{z}{\sqrt{\pi}}(1+\theta z^2)^pM(\frac{1}{2};\frac{3}{2};-z^2)
= \frac{z}{\sqrt{\pi}}\sum_{n=0}^{\infty} u_n (-z^2)^{n},
\end{equation*}
with $u_0 = 1$, $u_1 = \frac{1}{3} - \theta p$, $u_2 = \frac{1}{2} \left( \frac{1}{5} - \frac{2\theta p}{3} + \theta^2(p-1)p \right)$ and
\[
    u_{n+1} = \beta_0(n)u_n + \beta_1(n)u_{n-1} + \beta_2(n)u_{n-2},
\]
where
\begin{align*}
\beta_0(n) &= \frac{1 + 2\theta n(1+2n) - \theta p(3+4n) + n}{(n+1)(3+2n)}, \\
\beta_1(n) &= \frac{\theta(-1 - \theta(n-p-1)(\frac{3}{2}+n-p-2) - 2n + p + 2)}{(n+1)(\frac{3}{2}+n)}, \\
\beta_2(n) &= \frac{\theta^2(2n-2p-3)}{(n+1)(3+2n)},
\end{align*}
which could be written as:
\begin{equation*}
(1+\theta z^2)^p Erf(z) = \sum_{n=0}^{\infty}a_nz^{2n+1},
\end{equation*}
with $a_n = \frac{(-1)^nu_n}{\sqrt{\pi}}$, which completes the proof.

\end{proof}
\end{corollary}

Next, we present recurrence relations involving the lower(upper) incomplete gamma function $\gamma(s,z)$ ($\Gamma(s,z)$) respectively.
\begin{corollary}
Let $s,z \in \mathbb{C}$, $-(s+1)\notin \mathbb{N\cup }\left\{ 0\right\}$, then we have:
\begin{equation*}
(1-\theta z)^p\gamma(s,z) = (1-\theta z)^p\frac{z^s}{s}M(1;s+1;-z) = \displaystyle\sum_{n=0}^{\infty}\frac{(-1)^{n+s}u_n}{s}z^{n+s}
\end{equation*}
with $z \in \mathbb{C}$,$\quad |z| < \frac{1}{|\theta|}$, $|\arg (1-\theta z)| < \pi\quad and \quad |\arg z|<\pi $, $u_0 = 1$, $u_1 = \frac{1}{s+1} - \theta p$, $u_2 = \frac{1}{2} \left( \frac{2}{(s+1)^2+s+1} - \frac{2\theta p}{s+1} + \theta^2(p-1)p \right)$ and
\[
    u_{n+1} = \beta_0(n)u_n + \beta_1(n)u_{n-1} + \beta_2(n)u_{n-2},
\]
where
\begin{align*}
\beta_0(n) &= \frac{1 + 2\theta n(s+n) - \theta p(s+1+2n) + n}{(n+1)(s+1+n)}, \\
\beta_1(n) &= \frac{\theta(-2 - \theta(n-p-1)(s+n-p-1) - 2n + p + 2)}{(n+1)(s+1+n)}, \\
\beta_2(n) &= \frac{\theta^2(n-p-1)}{(n+1)(s+1+n)}.
\end{align*}

\begin{proof}
Applying Theorem \ref{thm2.2}, we just get:
\begin{equation*}
(1-\theta z)^pM(1;s+1;-z) = \displaystyle\sum_{n=0}^{\infty}u_n(-z)^n,
\end{equation*}
with $u_0 = 1$, $u_1 = \frac{1}{s+1} - \theta p$, $u_2 = \frac{1}{2} \left( \frac{2}{(s+1)^2+s+1} - \frac{2\theta p}{s+1} + \theta^2(p-1)p \right)$ and
\[
    u_{n+1} = \beta_0(n)u_n + \beta_1(n)u_{n-1} + \beta_2(n)u_{n-2},
\]
where
\begin{align*}
\beta_0(n) &= \frac{1 + 2\theta n(s+n) - \theta p(s+1+2n) + n}{(n+1)(s+1+n)}, \\
\beta_1(n) &= \frac{\theta(-2 - \theta(n-p-1)(s+n-p-1) - 2n + p + 2)}{(n+1)(s+1+n)}, \\
\beta_2(n) &= \frac{\theta^2(n-p-1)}{(n+1)(s+1+n)}.
\end{align*}
Then,
\begin{equation*}
(1-\theta z)^p\gamma(s,z) = \displaystyle\sum_{n=0}^{\infty}\frac{u_n}{s}(-z)^{n+s}
\end{equation*}
which completes the proof.
\end{proof}

\end{corollary}

\begin{corollary}
Let $s,z,p \in \mathbb{C}$, $\quad |z| < \frac{1}{|\theta|}$ with $\arg(1-\theta z) \quad and \quad \arg z \in (-\pi,\pi)$, then we have:
\begin{equation*}
(1- \theta z)^p\Gamma(s,z) = \sum_{n=0}^{\infty}v_nz^n + \sum_{n=0}^{\infty}w_nz^{n+s},
\end{equation*}
with $v_n = \Gamma(s)\sum_{k=0}^{n}\frac{(-1)^ku_{n-k}}{k!}$, $w_n = \frac{-1}{s}\sum_{k=0}^{n}\frac{(-1)^ku'_{n-k}}{k!}$, $u_0 = 1$, $u_1 = 1 - \theta p$, $u_2 = \frac{1}{2} \left( 1 - 2\theta p + \theta^2(p-1)p \right)$, $u'_0 = 1$, $u'_1 = \frac{1}{s+1} - \theta p$, $u'_2 = \frac{1}{2} \left( \frac{2}{(s+1)(s+2)} - \frac{2\theta p}{s+1} + \theta^2(p-1)p \right)$ and
\[
    u_{n+1} = \beta_0(n)u_n + \beta_1(n)u_{n-1} + \beta_2(n)u_{n-2},
\]
\[
    u'_{n+1} = \beta'_0(n)u'_n + \beta'_1(n)u'_{n-1} + \beta'_2(n)u'_{n-2},
\]
where
\begin{align*}
\beta_0(n) &= \frac{1-s + 2\theta n(n-s) - \theta p(1-s+2n) + n}{(n+1)(1-s+n)}, \\
\beta_1(n) &= \frac{\theta(-2(1-s) - \theta(n-p-1)(n-s-p-1) - 2n + p + 2)}{(n+1)(1-s+n)}, \\
\beta_2(n) &= \frac{\theta^2(n-s-p-1)}{(n+1)(1-s+n)}.
\end{align*}
\begin{align*}
\beta'_0(n) &= \frac{1 + 2\theta n(c+n-1) - \theta p(s+1+2n) + n}{(n+1)(s+1+n)}, \\
\beta'_1(n) &= \frac{\theta(-2 - \theta(n-p-1)(s+n-p-1) - 2n + p + 2)}{(n+1)(s+1+n)}, \\
\beta'_2(n) &= \frac{\theta^2(1+n-p-2)}{(n+1)(s+1+n)}.
\end{align*}

\begin{proof}
It is known to all that:
\begin{equation}\label{upper gamma}
\Gamma(s,z) = e^{-z}U(1-s;1-s;z),
\end{equation}
while $U(a;b;z) = \frac{\Gamma(1-b)}{\Gamma(a+1-b)}M(a;b;z) + \frac{\Gamma(b-1)}{\Gamma(a)}z^{1-b}M(a+1-b;2-b;z)$ is a combination of two solution of the Kummer's equation.

Now, substitute the the above equation into (\ref{upper gamma}), we get:
\begin{equation*}
\begin{aligned}
& (1- \theta z)^p\Gamma(s,z) = e^{-z}(1- \theta z)^p
(\frac{\Gamma(s)}{\Gamma(1)}M1-s;1-s;z) + \frac{\Gamma(-s)}{\Gamma(1-s)}z^{s}M(1;s+1;z))\\
& = e^{-z}(\Gamma(s)\sum_{n=0}^{\infty}u_nz^n + \frac{-1}{s}z^{s}\sum_{n=0}^{\infty}u'_nz^n) = \sum_{n=0}^{\infty}\frac{(-1)^nz^n}{n!}(\Gamma(s)\sum_{n=0}^{\infty}u_nz^n + \frac{-1}{s}z^{s}\sum_{n=0}^{\infty}u'_nz^n)\\
& = \sum_{n=0}^{\infty}\Gamma(s)\sum_{k=0}^{n}\frac{(-1)^ku_{n-k}}{k!}z^n + z^s\sum_{n=0}^{\infty}\frac{-1}{s}\sum_{k=0}^{n}\frac{(-1)^ku'_{n-k}}{k!}z^{n}
\end{aligned}
\end{equation*}
while the explicit expression of $u_n$ and $u'_n$ could be calculated directly by Theorem \ref{thm2.2}, we could just get the final expression:
\begin{equation*}
\begin{aligned}
(1- \theta z)^p\Gamma(s,z)= \sum_{n=0}^{\infty}v_nz^n + \sum_{n=0}^{\infty}w_nz^{n+s}
\end{aligned}
\end{equation*}
with $v_n = \Gamma(s)\sum_{k=0}^{n}\frac{(-1)^ku_{n-k}}{k!}$, $w_n = \frac{-1}{s}\sum_{k=0}^{n}\frac{(-1)^ku'_{n-k}}{k!}$,
which completes the proof.

\end{proof}
\end{corollary}

We further state a corollary involving the modified Bessel function $I_{a-\frac{1}{2}}(z)$.
\begin{corollary}
Let $a,\theta, p, z \in \mathbb{C}$, $-2a\notin \mathbb{N\cup }\left\{ 0\right\}$,$\quad |z| < \frac{1}{|\theta|}$ and $|\arg (1-\theta z)| < \pi$,  then we have:
\begin{equation*}
(1-\theta z)^p I_{a-\frac{1}{2}}(z) = \displaystyle\sum_{n=0}^{\infty}a_nz^n
\end{equation*}
where
\begin{equation*}
    a_n = \frac{4^{\frac{1}{2}-a}}{\Gamma(a+\frac{1}{2})}\displaystyle\sum_{k=0}^{n}
    \frac{(-1)^k u_{n-k}}{2^k k!},
\end{equation*}
\begin{equation*}
\begin{aligned}
    & u_0 = 1, u_1 = \frac{1}{2} - \theta p, u_2 = \frac{1}{2} \left( \frac{a+1}{2a+1} - \theta p + \theta^2(p-1)p \right)\\
    & u_{n+1} = \beta_0(n)u_n + \beta_1(n)u_{n-1} + \beta_2(n)u_{n-2},\\
    & \beta_0(n) = \frac{a + 2\theta n(2a+n-1) - \theta p(2a+2n) + n}{(n+1)(2a+n)}, \\
    & \beta_1(n) = \frac{\theta(-2a - \theta(n-p-1)(2a+n-p-2) - 2n + p + 2)}{(n+1)(2a+n)}, \\
    & \beta_2(n) = \frac{\theta^2(a+n-p-2)}{(n+1)(2a+n)}.
\end{aligned}
\end{equation*}

\begin{proof}
It is known to all that:
\begin{equation*}
M(a;2a;z) = I_{a-\frac{1}{2}}(z)\Gamma(a+\frac{1}{2})
(\frac{z}{4})^{\frac{1}{2} - a}e^\frac{z}{2}
\end{equation*}
Therefore, applying Theorem \ref{thm2.2} gives
\begin{equation*}
\begin{aligned}
    & (1-\theta z)^p I_{a-\frac{1}{2}}(z) = e^{-\frac{z}{2}}
    \frac{4^{\frac{1}{2}-a}}{\Gamma(a+\frac{1}{2})}z^{a-\frac{1}{2}}
    (1-\theta z)^pM(a;2a;z)\\
    &  =  \frac{4^{\frac{1}{2}-a}}{\Gamma(a+\frac{1}{2})}z^{a-\frac{1}{2}} \displaystyle\sum_{n=0}^{\infty}\frac{(-1)^nz^n}{2^nn!}
    \displaystyle\sum_{n=0}^{\infty}u_nz^n\\
    &  = \frac{4^{\frac{1}{2}-a}}{\Gamma(a+\frac{1}{2})}z^{a-\frac{1}{2}} \displaystyle\sum_{n=0}^{\infty} \displaystyle\sum_{k=0}^{n}
    \frac{(-1)^k u_{n-k}}{2^k k!}z^n\\
    & = \displaystyle\sum_{n=0}^{\infty}a_nz^{n+a-\frac{1}{2}},
\end{aligned}
\end{equation*}
where $u_n$ could be calculated through Theorem \ref{thm2.2} and $a_n = \frac{4^{\frac{1}{2}-a}}{\Gamma(a+\frac{1}{2})}\displaystyle\sum_{k=0}^{n}
    \frac{(-1)^k u_{n-k}}{2^k k!}$, which completes the proof.
\end{proof}

\end{corollary}

%------corollary for gauss--------
\subsection{Applications of the Gauss hypergeometric function}
\mbox{} % 这里加一个空盒子并强制换行

First, we present the recurrence relation involving zero-balanced hypergeometric function by substituting $c=a+b$ into Theorem \ref{main thm1}.
\begin{corollary}
Let $a,b,z \in \mathbb{C}$, $-(a+b))\notin \mathbb{N\cup }\left\{ 0\right\} $. Then we have:
\begin{equation*}
(1-\theta z)^pF(a,b;a+b;z) = \sum_{n=0}^{\infty }v_{n}z^{n}
\end{equation*}
with $|z| < 1$, $|\arg (1-\theta z)| < \pi$ and $v_{0}=1$, $v_{1}=\frac{ab}{a+b}-p\theta $,
\begin{equation*}
v_{2}=\frac{1}{2}\theta ^{2}p\left( p-1\right) -\theta p\frac{ab}{a+b}+\frac{1
}{2}\frac{ab\left( b+1\right) \left( a+1\right) }{(a+b)\left( a+b+1\right) },
\end{equation*}%
and for $n\geq 2$,%
\begin{equation*}
v_{n+1}=\frac{\mu _{n,p,\theta }\left( a,b\right) }{\left( n+1\right)
\left( n+a+b\right) }v_{n}-\theta \dfrac{\zeta _{n,p,\theta }\left(
a,b\right) }{\left( n+1\right) \left( n+a+b\right) }v_{n-1}+\theta ^{2}\frac{
\sigma _{n,p}\left( a,b\right) }{\left( n+1\right) \left( n+a+b\right) }
v_{n-2},
\end{equation*}
where
\begin{eqnarray*}
\mu _{n,p,\theta }\left( a,b\right) &=&\left( n+a\right) \left( n+b\right)
+\theta \left[ 2n^{2}-2n\left( p-(a+b)+1\right) -(a+b)p\right] , \\
\zeta _{n,p,\theta }\left( a,b\right) &=&2n^{2}+2\left( a+b-p-2\right)
n-\left( a+b-1\right) p+2\left( a-1\right) \left( b-1\right) \\
&&+\theta \left( n-p-1\right) \left( n-p+a+b-2\right) , \\
\sigma _{n,p}\left( a,b\right) &=&\left( n+a-p-2\right) \left(
n+b-p-2\right).
\end{eqnarray*}

\end{corollary}

Next, we present a recurrence relation involving the product of the logarithmic function and the Gauss hypergeometric function.
\begin{corollary}
\label{C-U0}Let $a,b,c,z \in \mathbb{C}$, $-c\notin \mathbb{N\cup }\left\{ 0\right\} $. Then we have
\begin{equation*}
V\left( z\right) =\ln \left( 1- z\right) \times F\left( a,b;c;z\right)
=\sum_{n=0}^{\infty }\lambda_{n}z^{n},  \label{V}
\end{equation*}
with $|z| < 1$, $|\arg(1- z)| < \pi$, $\lambda_{0}=0$, $\lambda_{1}=-1$ and for $n\geq 1$,%
\begin{equation}
\lambda_{n+1}=\tau _{n}\lambda_{n}-\rho_{n}\lambda_{n-1}+\gamma _{n}w_{n},
\label{vn+1-n-1-rr}
\end{equation}%
where $\tau _{n}$, $\rho _{n}$ are defined in (\ref{aln}), (\ref{ben}),
and
\begin{eqnarray*}
\gamma _{n} &=&\frac{\left( c-b-a\right) n^{2}+\left( a+b-2ab\right)
n-c\left( a-1\right) \left( b-1\right) }{\left( n+1\right) \left(
n+a-1\right) \left( n+b-1\right) \left( n+c\right) }, \\
w_{n} &=&\frac{\left( a\right) _{n}\left( b\right) _{n}}{n!\left( c\right)
_{n}}=\frac{\Gamma \left( c\right) }{\Gamma \left( a\right) \Gamma \left(
b\right) }\frac{\Gamma \left( n+a\right) \Gamma \left( n+b\right) }{n!\Gamma
\left( n+c\right) }.
\end{eqnarray*}
\end{corollary}

\begin{proof}
Notice that
\begin{equation*}
F\left( a,b;c;z\right) =\sum_{n=0}^{\infty }w_{n}z^{n}\text{ \ and \ }%
\lim_{p\rightarrow 0}\frac{\left( 1- z\right) ^{p}-1}{p}=\ln \left(1- z\right),
\end{equation*}
we have
\begin{equation*}
V\left( z\right) =\ln \left( 1- z\right) \times F\left( a,b;c;z\right)
=\lim_{p\rightarrow 0}\frac{\left( 1- z\right) ^{p}F\left( a,b;c;z\right)
-F\left( a,b;c;z\right) }{p},
\end{equation*}%
and then,
\begin{equation*}
\sum_{n=0}^{\infty }\lambda_{n}z^{n}=\lim_{p\rightarrow 0}\sum_{n=0}^{\infty }%
\frac{v_{n}-w_{n}}{p}z^{n}=\sum_{n=0}^{\infty }\left( \lim_{p\rightarrow 0}%
\frac{v_{n}-w_{n}}{p}\right) z^{n},
\end{equation*}%
which implies that%
\begin{equation*}
\lambda_{n}=\lim_{p\rightarrow 0}\frac{v_{n}-w_{n}}{p}\text{ \ for }n\geq 0\text{.}
\end{equation*}%
Evidently, $\lambda_{0}=0$, $\lambda_{1}=-1$. To obtain (\ref{vn+1-n-1-rr}), we rewrite
the recurrence relation (\ref{u-th=1-rr}) as%
\begin{equation*}
\begin{array}{c}
\dfrac{v_{n+1}-w_{n+1}}{p}=\tau _{n}\dfrac{v_{n}-w_{n}}{p}-\rho _{n}%
\dfrac{v_{n-1}-w_{n-1}}{p}
+\dfrac{\tau _{n}w_{n}-\rho _{n}w_{n-1}-w_{n+1}}{p}.
\end{array}
\label{un-wn/p-rr}
\end{equation*}%
Noting that%
\begin{equation*}
\frac{\tau _{n}w_{n}-\rho _{n}w_{n-1}-w_{n+1}}{p}=\frac{w_{n}}{p}\left[
\tau _{n}-\rho _{n}\frac{n\left( n+c-1\right) }{\left( n+a-1\right)
\left( n+b-1\right) }-\frac{\left( n+a\right) \left( n+b\right) }{\left(
n+1\right) \left( n+c\right) }\right]
\end{equation*}%
\begin{eqnarray*}
&=&\frac{\left( c-b-a-p\right) n^{2}+\left( a+b-2ab+p-cp\right) n-c\left(
a-1\right) \left( b-1\right) }{\left( n+1\right) \left( n+a-1\right) \left(
n+b-1\right) \left( n+c\right) }w_{n} \\
&\rightarrow &\frac{\left( c-b-a\right) n^{2}+\left( a+b-2ab\right)
n-c\left( a-1\right) \left( b-1\right) }{\left( n+1\right) \left(
n+a-1\right) \left( n+b-1\right) \left( n+c\right) }w_{n}=\gamma _{n}w_{n}%
\text{,}
\end{eqnarray*}%
the required recurrence relation (\ref{vn+1-n-1-rr}) just follows as $p\rightarrow 0$, which completes the proof.
\end{proof}

We also obtain a special case corresponding to the complete elliptic integrals of the first and second kinds by substituting the explicit $a,b,c$ into Theorem\ref{main thm1}.
\begin{corollary}
Let $0<r<1$, then we have:
\begin{equation*}
    (1 - \theta r^{2})^{p}K(r) = \frac{\pi}{2}
     (1 - \theta r^{2})^{p}F\!\left( \frac{1}{2}, \frac{1}{2}; 1; r^2 \right) = \displaystyle\sum_{n=0}^{\infty}
     \frac{\pi}{2}v_n r^{2n}
\end{equation*}
\begin{equation*}
        (1 - \theta r^{2})^{p}E(r) = \frac{\pi}{2}
     (1 - \theta r^{2})^{p}F\!\left( -\frac{1}{2}, \frac{1}{2}; 1; r^2 \right) = \displaystyle\sum_{n=0}^{\infty}
     \frac{\pi}{2}v'_n r^{2n},
\end{equation*}
with $v_0 = 1$, $v_1 = \frac{1}{4} - p \theta$,
$v_2 = \frac{1}{2} \theta ^{2}p(p-1) - \frac{1}{4} \theta p
+ \frac{9}{64}$ and $v'_0 = 1$, $v'_1 = -\frac{1}{4} - p \theta$,
$v'_2 = \frac{1}{2} \theta ^{2}p(p-1) + \frac{1}{4} \theta p
- \frac{3}{64}$,

and for $n \ge 2$,
\begin{equation*}
v_{n+1} = \frac{\mu \left( n,p,\theta \right) }{\left( n + 1\right)^2}v_{n}
- \theta \dfrac{\zeta \left( n,p,\theta \right) }{\left( n+1\right)^2}v_{n-1} +
\theta ^{2}\frac{\sigma \left( n,p \right) }{\left( n+1\right)^2}
v_{n-2},
\end{equation*}
\begin{equation*}
v'_{n+1} = \frac{\mu' \left( n,p,\theta \right) }{\left( n + 1\right)^2}v'_{n}
- \theta \dfrac{\zeta' \left( n,p,\theta \right) }{\left( n+1\right)^2}v'_{n-1} +
\theta ^{2}\frac{\sigma' \left( n,p \right) }{\left( n+1\right)^2}v'_{n-2},
\end{equation*}

where
\begin{eqnarray*}
\mu \left( n,p,\theta \right) &=&\left( n+\frac{1}{2}\right) \left( n+\frac{1}{2}\right)
+\theta \left[ 2n^{2}-2np -p\right] , \\
\zeta \left( n,p,\theta \right) &=&2n^{2}-2\left( p+1\right)
n + \frac{1}{2}+\theta \left( n-p-1\right)^2, \\
\sigma \left( n,p \right) &=&\left( n-p-\frac{3}{2}\right)^2,
\end{eqnarray*}

\begin{eqnarray*}
\mu' \left( n,p,\theta \right) &=& n^2 - \frac{1}{4} + \theta \left( 2n^{2}-2np - p\right) , \\
\zeta' \left( n,p,\theta \right) &=&2n^{2}-2\left( p + 2 \right)
n + p - \frac{3}{2} + \theta \left( n-p-1\right)^2, \\
\sigma' \left( n,p \right) &=&\left( n-p-\frac{5}{2}\right) \left(n-p-\frac{3}{2}\right).
\end{eqnarray*}
\end{corollary}

We further provide an explicit power series expansion for the product of a logarithmic function and a power function.
\begin{corollary}
Let $z,p \in \mathbb{C}$, $-p \notin \mathbb{N}$, $|z| < 1$ and $|\arg(1-z)|,|\arg(1+z)| < \pi $, then we have
\begin{equation*}
(1+z)^{p} ln(1-z) = - \displaystyle\sum_{n=0}^{\infty}
v_nz^{n+1}
\end{equation*}
where
\begin{equation}\label{equ3.2}
v_n = \frac{(-1)^n (1-p)_n}{(n+1)!}
\sum_{k=0}^{n} \frac{(-1)^k k!}{(1-p)_k}
\left[ \sum_{j=0}^{k} \frac{(-1)^j (-p-1)_j}{j!} \right].
\end{equation}
If there is a $-p_0 \in \mathbb{N}$ such that $(1-p_0)_k = 0$ for certain $0 \le k \le n$, then we define $v_n = \lim_{p \to p_0} v_n(p)$.
\begin{proof}
We only consider the case of $-p \notin \mathbb{N}$. Note that:
\begin{equation*}
-\frac{\ln(1-z)}{z}
= \sum_{n=0}^{\infty} \frac{z^n}{n+1}= F(1,1;2;z).
\end{equation*}
Taking $(a,b,c) = (1,1,2)$ in Corollary \ref{col2.1} gives $v_0 = 1$, $v_1 = \tfrac{1}{2} + p$, $v_2 = \tfrac{3p^2 + 2}{6}$ and for $n \ge 2$,
\begin{equation*}
v_{n+1} = -\frac{n - 2p - 1}{n+2} v_n + \frac{n^2 + n - p^2 - 2p}{(n+1)(n+2)} v_{n-1} + \frac{(n-p-1)^2}{(n+1)(n+2)} v_{n-2}.
\end{equation*}
This 3-order recurrence relation can be reduced to a 2-order one:
\begin{equation}\label{equ3.5}
 v'_{n+1} = \frac{p+2}{n+2} v'_n + \frac{n(n-p-1)}{(n+2)(n+1)} v'_{n-1}, \qquad n \ge 1.
\end{equation}
where
\begin{equation}\label{equ3.6}
v'_n = v_n + \frac{n - p}{n + 1} v_{n-1}
\end{equation}
with $v'_1 = (p+2)/2$ and $v'_2 = (p^2 + 3p + 4)/6$.
Furthermore, the 2-order recurrence relation (\ref{equ3.5}) can be reduced to 1-order one:
\begin{equation}\label{equ3.7}
v''_{n+1} = -\frac{n - p - 1}{n + 2} v''_n,
\qquad n \ge 2,
\end{equation}
where
\begin{center}
$v''_n = v'_n - \frac{n}{n + 1} v'_{n-1}$ with
$v''_2 = \frac{1}{6} p (p + 1)$.
\end{center}
Solving the recurrence equation (\ref{equ3.7}) gives
\[
     v'_n - \frac{n}{n + 1} v'_{n-1} = v''_n
     = \frac{(-1)^n (-1 - p)_n}{(n + 1)!},
\]
which, in combination with (\ref{equ3.6}), implies that
\[
    (k + 1) v_k + (k - p) v_{k-1} = (k + 1) v'_k
    = \sum_{j=0}^{k} \frac{(-1)^j (-1 - p)_j}{(j + 1)!},
\]
for $k \geq 1$. Solving this for $v_k$ leads to (\ref{equ3.2}), which completes the proof.
\end{proof}

\end{corollary}

The following corollary provides a recurrence relation for the coefficients of the product of $(1-z^2)$ and the inverse hyperbolic tangent function $arctanhz$.
\begin{corollary}
Let $z \in \mathbb{C}$, $|z| <1$, then we have:
\begin{equation*}
f(z) = (1-z^2)arctanhz =
\displaystyle\sum_{n=0}^{\infty} v_nz^{2n+1}
\end{equation*}
with $v_0 = 1$, $v_1 = -\frac{1}{6}$, and for $n \geq 1$,
\begin{equation*}
v_n = \tau_nv_n - \rho_nv_{n-1},
\end{equation*}
where
\begin{equation*}
\tau_n = \frac{2n^2-1}{(n + \frac{3}{2})(n+1)}
\end{equation*}
\begin{equation*}
\rho_n = \frac{(n-1)(2n-3)}{(n+1)(2n+3)}
\end{equation*}

\begin{proof}
It is know to all that we have a Gauss hypergeometric representation for the inverse hyperbolic tangent function:
\begin{equation*}
arctanhz = z F(\frac{1}{2},1;\frac{3}{2};z^2)
\end{equation*}
Therefore, by Corollary \ref{col2.3}, the multiplication becomes:
\begin{equation*}
f(z) = z(1-z^2)F(\frac{1}{2},1;\frac{3}{2};z^2) =
z\displaystyle\sum_{n=0}^{\infty}v_nz^{2n}.
\end{equation*}
And the coefficient $v_n$ can be calculated through Corollary \ref{col2.3}, which completes the proof.
\end{proof}
\end{corollary}

The following corollary provides a recurrence relation for the coefficients of Chebyshev polynomials $y=T_m(z) = T_n(\cos\theta) = \cos(n \theta),z \in \mathbb{C}$. Chebyshev polynomials are the solutions to the following equation:
\begin{equation*}
(1-z^2)\frac{d^2y}{dz^2} - z \frac{dy}{dz} + m^2y = 0
\end{equation*}
And it has a representation in Gauss hypergeometric function:
\begin{equation}\label{cheby}
T_m(z) = F(-m,m;\frac{1}{2};\frac{1-z}{2})
\end{equation}
\begin{corollary}
Let $z,\theta,p \in \mathbb{C}$, $|\frac{1-z}{2}|<1$ and $|\arg (1-\theta\frac{1-z}{2})|< \pi$, then we have:
\begin{equation*}
f(z) = (1-\theta\frac{1-z}{2})^pT_m(z) = \displaystyle\sum_{n=0}^{\infty} v_n (\frac{1-z}{2})^n
\end{equation*}
with $v_{0}=1$, $v_{1}=-2m^2-p\theta $,
\begin{equation*}
v_{2}=\frac{1}{2}\theta ^{2}p\left( p-1\right) + 2m^2\theta p+\frac{4m^2(m+1)(m-1)}{6},
\end{equation*}
and for $n\geq 2$,
\begin{equation*}
v_{n+1}=\frac{\mu _{n,p,\theta }\left(m\right) }{\left( n+1\right)
\left( n+\frac{1}{2}\right) }v_{n}-\theta \dfrac{\zeta _{n,p,\theta }\left(
m\right) }{\left( n+1\right) \left( n+\frac{1}{2}\right) }v_{n-1}+\theta ^{2}\frac{%
\sigma _{n,p}\left( m\right) }{\left( n+1\right) \left( n+\frac{1}{2}\right) }%
v_{n-2},
\end{equation*}%
where
\begin{eqnarray*}
\mu _{n,p,\theta }\left( m\right) &=&\left( n-m\right) \left( n+m\right)
+\theta \left[ 2n^{2}-2n\left( p-c+1\right) -cp\right] , \\
\zeta _{n,p,\theta }\left( m\right) &=&2n^{2}-2\left(p+2\right)
n + p+2\left( -m-1\right) \left( m-1\right) \\
&&+\theta \left( n-p-1\right) \left( n-p+\frac{1}{2}-2\right) , \\
\sigma _{n,p}\left( m\right) &=&\left( n-m-p-2\right) \left(
n+m-p-2\right).
\end{eqnarray*}
\begin{proof}
We just use equation (\ref{cheby}) to derive:
\begin{equation*}
(1-\theta\frac{1-z}{2})^pT_m(z) = (1-\theta\frac{1-z}{2})^pF(-m,m;\frac{1}{2};\frac{1-z}{2}).
\end{equation*}
Then we just apply Theorem \ref{main thm1} to calculate the explicit expression of recurrence relations, which completes the proof.
\end{proof}
\end{corollary}

The following corollary provides a recurrence relation for the coefficients of the product of $(1+\theta z)^p$ and $(1+z)^{\alpha}$.
\begin{corollary}
Let $z,\theta,p,\alpha,\beta \in \mathbb{C}$, $|z|<1$ and $|\arg (1+\theta z)|<\pi$, then we have:
\begin{equation*}
f(z) = (1+\theta z)^p(1+z)^{\alpha} =
\displaystyle\sum_{n=0}^{\infty}v_n(-z)^n
\end{equation*}
with $v_{0}=1$, $v_{1}=\alpha-p\theta $,%
\begin{equation*}
v_{2}=\frac{1}{2}\theta ^{2}p\left( p-1\right) -\theta p\alpha+\alpha(\alpha+1),
\end{equation*}
and for $n\geq 2$,
\begin{equation}
v_{n+1}=\frac{\mu _{n,p,\theta }\left( \alpha,\beta\right) }{\left( n+1\right)
\left( n+\beta\right) }v_{n}-\theta \dfrac{\zeta _{n,p,\theta }\left(
\alpha,\beta\right) }{\left( n+1\right) \left( n+\beta\right) }v_{n-1}+\theta ^{2}\frac{%
\sigma _{n,p}\left( \alpha,\beta\right) }{\left( n+1\right) \left( n+\beta\right) }%
v_{n-2},
\end{equation}%
where
\begin{eqnarray*}
\mu _{n,p,\theta }\left( \alpha,\beta\right) &=&\left( n+\alpha\right) \left( n+\beta\right)
+\theta \left[ 2n^{2}-2n\left( p-\beta+1\right) -\beta p\right] , \\
\zeta _{n,p,\theta }\left( \alpha,\beta\right) &=&2n^{2}+2\left( \alpha+\beta-p-2\right)
n-\left( \alpha+\beta-1\right) p+2\left( \alpha-1\right) \left( \beta-1\right) \\
&&+\theta \left( n-p-1\right) \left( n-p+\beta-2\right) , \\
\sigma _{n,p}\left( \alpha,\beta\right) &=&\left( n+\alpha-p-2\right) \left(
n+\beta-p-2\right).
\end{eqnarray*}
\begin{proof}
It is know to all that $(1+z)^{\alpha} = F(\alpha,\beta;\beta;-z)$, then we substitute it to the above equation and apply Theorem \ref{main thm1}, which completes the proof.
\end{proof}
\end{corollary}

The following corollary provides a recurrence relation for the coefficients of the product of $(1-\theta z^2)^p$ and $\arcsin z$.
\begin{corollary}
Let $z,\theta,p \in \mathbb{C}$, $|z|<1$ and $|\arg (1-\theta z^2)|<\pi$, then we have:
\begin{equation*}
f(z) = (1-\theta z^2)^p\arcsin z = \displaystyle\sum_{n=0}^{\infty}v_nz^{2n+1}
\end{equation*}
with $v_{0}=1$, $v_{1}=\frac{1}{6}-p\theta $,%
\begin{equation*}
v_{2}=\frac{1}{2}\theta ^{2}p\left( p-1\right) -\frac{1}{6}\theta p+
\frac{3}{40},
\end{equation*}
and for $n\geq 2$,
\begin{equation}
v_{n+1}=\frac{\mu _{n,p,\theta } }{\left( n+1\right)
\left( n+\frac{3}{2}\right) }v_{n}-\theta \dfrac{\zeta _{n,p,\theta }}{\left( n+1\right) \left( n+\frac{3}{2}\right) }v_{n-1}+\theta ^{2}\frac{%
\sigma _{n,p}}{\left( n+1\right) \left( n+\frac{3}{2}\right) }%
v_{n-2},
\end{equation}%
where
\begin{eqnarray*}
\mu _{n,p,\theta } &=&\left( n+\frac{1}{2}\right) \left( n+\frac{1}{2}\right)
+\theta \left[ 2n^{2}-2n\left( p-\frac{1}{2}\right) -\frac{3}{2}p\right] , \\
\zeta _{n,p,\theta } &=&2n^{2}-2\left(p+1\right)
n+\frac{1}{2} +\theta \left( n-p-1\right) \left( n-p-\frac{1}{2}\right) , \\
\sigma _{n,p} &=&\left( n-p-\frac{1}{2}\right) \left(
n-p-\frac{1}{2}\right).
\end{eqnarray*}
\begin{proof}
It is known to all that $\arcsin z = zF(\frac{1}{2},\frac{1}{2};\frac{3}{2};z^2)$, then we substitute it into the above equation and apply Theorem \ref{main thm1}, which completes the proof.
\end{proof}
\end{corollary}

The following corollary provides a recurrence relation for the coefficients of the product of $(1+\theta z^2)^p$ and $\arctan z$.
\begin{corollary}
Let $z,\theta,p \in \mathbb{C}$, $|z|<1$ and $|\arg (1+\theta z^2)|<\pi$, then we have:
\begin{equation*}
f(z) = (1+\theta z^2)^p\arctan z = \displaystyle\sum_{n=0}^{\infty}(-1)^nv_nz^{2n+1}
\end{equation*}
with $v_{0}=1$, $v_{1}=\frac{1}{3}-p\theta $,%
\begin{equation*}
v_{2}=\frac{1}{2}\theta ^{2}p\left( p-1\right) -\frac{1}{3}\theta p+\frac{1}{5},
\end{equation*}
and for $n\geq 2$,
\begin{equation}
v_{n+1}=\frac{\mu _{n,p,\theta } }{\left( n+1\right)
\left( n+\frac{3}{2}\right) }v_{n}-\theta \dfrac{\zeta _{n,p,\theta } }{\left( n+1\right) \left( n+\frac{3}{2}\right) }v_{n-1}+\theta ^{2}\frac{%
\sigma _{n,p} }{\left( n+1\right) \left( n+\frac{3}{2}\right) }%
v_{n-2},
\end{equation}%
where
\begin{eqnarray*}
\mu _{n,p,\theta } &=&\left( n+\frac{1}{2}\right) \left( n+1\right)
+\theta \left[ 2n^{2}-2n\left( p-\frac{1}{2}\right) -\frac{3}{2}p\right] , \\
\zeta _{n,p,\theta } &=&2n^{2}-2\left( p+\frac{1}{2}\right)
n - \frac{1}{2}p + \theta \left( n-p-1\right) \left( n-p-\frac{1}{2}\right) , \\
\sigma _{n,p} &=&\left( n+\frac{1}{2}-p-2\right) \left(
n-p-1\right).
\end{eqnarray*}
\begin{proof}
It is known to all that $\arctan z = zF(\frac{1}{2},1;\frac{3}{2};-z^2)$, then we substitute it into the above equation and apply Theorem \ref{main thm1}, which completes the proof.
\end{proof}
\end{corollary}

The following corollary provides a recurrence relation for the coefficients of the product of $(1-\theta \sin^2 z)^p$ and $\cos \xi z$.
\begin{corollary}
Let $z,\xi,\theta,p \in \mathbb{C}$, $|\sin z|<1$ and $|\arg (1-\theta \sin^2 z)|<\pi$, then we have:
\begin{equation*}
f(z) = (1-\theta \sin^2 z)^p\cos \xi z = \displaystyle\sum_{n=0}^{\infty}v_n(\sin z)^{2n}
\end{equation*}
with $v_{0}=1$, $v_{1}=-\frac{\xi^2}{2}-p\theta $,%
\begin{equation*}
v_{2}=\frac{1}{2}\theta ^{2}p\left( p-1\right) +\frac{1}{2}p\theta \xi^2 +\frac{\xi^2\left( \xi^2-4\right) }{24},
\end{equation*}
and for $n\geq 2$,
\begin{equation*}
v_{n+1}=\frac{\mu _{n,p,\theta }}{\left( n+1\right)
\left( n+\frac{1}{2}\right) }v_{n}-\theta \dfrac{\zeta _{n,p,\theta }}{\left( n+1\right) \left( n+\frac{1}{2}\right) }v_{n-1}+\theta ^{2}\frac{%
\sigma _{n,p}}{\left( n+1\right) \left( n+\frac{1}{2}\right) }%
v_{n-2},
\end{equation*}%
where
\begin{eqnarray*}
\mu _{n,p,\theta } &=&n^2-\frac{\xi^2}{4}
+\theta \left[ 2n^{2}-\left( 2p+1\right) n -\frac{p}{2}\right] , \\
\zeta _{n,p,\theta } &=&2n^{2}-2\left( p+2\right)
n+p+2-\frac{\xi^2}{2} +\theta \left( n-p-1\right) \left( n-p-\frac{3}{2}\right) , \\
\sigma _{n,p} &=&\left( n-p-2\right)^2 -\frac{\xi^2}{4}.
\end{eqnarray*}
\begin{proof}
It is known to all that $\cos \mu z = F(\frac{\xi}{2},-\frac{\xi}{2};\frac{1}{2};\sin^2 z)$, then we substitute it into the above equation and apply Theorem \ref{main thm1}, which completes the proof.
\end{proof}
\end{corollary}

The following corollary provides a recurrence relation for the coefficients of the product of $(1-\theta \sin^2 z)^p$ and $\sin \xi z$.
\begin{corollary}
Let $z,\xi,\theta,p \in \mathbb{C}$, $|\sin z|<1$ and $|\arg (1-\theta \sin^2 z)|<\pi$, then we have:
\begin{equation*}
f(z) = (1-\theta \sin^2 z)^p\sin \xi z = \displaystyle\sum_{n=0}^{\infty}v_n(\sin z)^{2n+1}
\end{equation*}
with $v_{0}=1$, $v_{1}=\frac{1-\xi^2}{6}-p\theta $,%
\begin{equation*}
v_{2}=\frac{1}{2}\theta ^{2}p\left( p-1\right) +\frac{1}{6}p\theta \left( \xi^2-1\right) +\frac{\left( 1-\xi^2\right) \left( 9-\xi^2\right) }{120},
\end{equation*}
and for $n\geq 2$,
\begin{equation}
v_{n+1}=\frac{\mu _{n,p,\theta }}{\left( n+1\right)
\left( n+\frac{3}{2}\right) }v_{n}-\theta \dfrac{\zeta _{n,p,\theta }}{\left( n+1\right) \left( n+\frac{3}{2}\right) }v_{n-1}+\theta ^{2}\frac{%
\sigma _{n,p}}{\left( n+1\right) \left( n+\frac{3}{2}\right) }%
v_{n-2},
\end{equation}%
where
\begin{eqnarray*}
\mu _{n,p,\theta } &=&n^2+n+\frac{1-\xi^2}{4}
+\theta \left[ 2n^{2}-\left( 2p-1\right) n -\frac{3p}{2}\right] , \\
\zeta _{n,p,\theta } &=&2n^{2}-2\left( p+1\right)
n+\frac{1-\xi^2}{2} \\
&&+\theta \left( n-p-1\right) \left( n-p-\frac{1}{2}\right) , \\
\sigma _{n,p} &=&\left( n-p-\frac{3}{2}\right)^2 -\frac{\xi^2}{4}.
\end{eqnarray*}
\begin{proof}
It is known to all that $\sin \mu z = \mu \sin zF(\frac{1+\xi}{2},\frac{1-\xi}{2};\frac{3}{2};\sin^2 z)$, then we substitute it into the above equation and apply Theorem \ref{main thm1}, which completes the proof.
\end{proof}
\end{corollary}

The following corollary provides a recurrence relation for the coefficients of the product of $(1-\theta z^2)^p$ and $\log\frac{1+z}{1-z}$.
\begin{corollary}
Let $z,,\theta,p \in \mathbb{C}$, $|z|<1$ and $|\arg (1-\theta z^2)|,|\arg (\frac{1+z}{1-z})|<\pi$, then we have:
\begin{equation*}
f(z) = (1-\theta z^2)^p\log\frac{1+z}{1-z} = \displaystyle\sum_{n=0}^{\infty}2v_nz^{2n+1}
\end{equation*}
with $v_{0}=1$, $v_{1}=\frac{1}{3}-p\theta $,%
\begin{equation*}
v_{2}=\frac{1}{2}\theta ^{2}p\left( p-1\right) -\frac{1}{3}p\theta +\frac{1}{5},
\end{equation*}
and for $n\geq 2$,
\begin{equation}
v_{n+1}=\frac{\mu _{n,p,\theta }}{\left( n+1\right)
\left( n+\frac{3}{2}\right) }v_{n}-\theta \dfrac{\zeta _{n,p,\theta }}{\left( n+1\right) \left( n+\frac{3}{2}\right) }v_{n-1}+\theta ^{2}\frac{%
\sigma _{n,p}}{\left( n+1\right) \left( n+\frac{3}{2}\right) }%
v_{n-2},
\end{equation}%
where
\begin{eqnarray*}
\mu _{n,p,\theta } &=&n^2+\frac{3}{2}n+\frac{1}{2}
+\theta \left[ 2n^{2}-\left( 2p-1\right) n -\frac{3p}{2}\right] , \\
\zeta _{n,p,\theta } &=&2n^{2}-\left( 2p+1\right)
n-\frac{p}{2} \\
&&+\theta \left( n-p-1\right) \left( n-p-\frac{1}{2}\right) , \\
\sigma _{n,p} &=&\left( n-p-1\right)\left( n-p-\frac{3}{2}\right).
\end{eqnarray*}
\begin{proof}
It is known to all that $\log\frac{1+z}{1-z} = 2zF(\frac{1}{2},1;\frac{3}{2};z^2)$, then we just substitute it into the above equation and apply Theorem \ref{main thm1}, which completes the proof.
\end{proof}
\end{corollary}

Finally, we provide a recurrence relation of coefficients for the product of the inverse trigonometric function and inverse hyperbolic function.
\begin{corollary}
Let $z \in \mathbb{C}$, $|z|<1$ and $|\arg (1-z^2)|<\pi$, then we have:
\begin{equation*}
    f(z) = arcsinz \times arctanhz = \displaystyle\sum_{n=0}^{\infty}v_n z^{2n+2}
\end{equation*}
with $v_n = \frac{v_n'+v_n''}{2n+2}$ and $v_0' = v_0'' = 1$, $v_1' = \frac{5}{6}$, $v_1'' = \frac{7}{6}$, where
\begin{equation*}
\begin{aligned}
        &v_{n+1}' = \frac{2n^2+3n+\frac{5}{4}}{(n+1)(n+\frac{3}{2})}v_n' - \frac{n(2n+1)}{(n+1)(2n+3)}v_{n-1}'\\
        &u_{n+1}'' = \frac{2n^2+\frac{7}{2}n+\frac{7}{4}}{(n+1)(n+\frac{3}{2})}v_n'' - \frac{(n+\frac{1}{2})^{2}}{(n+1)(n+\frac{3}{2})}v_{n-1}''
    \end{aligned}
\end{equation*}

\begin{proof}
For the $f(z) = arcsinz \times arctanhz$, we can first differentiate it and get:
\begin{equation}\label{equ3.33}
\begin{aligned}
    &f' = (1-z^2)^{-\frac{1}{2}}arctanhz+(1-z^2)^{-1}arcsinz \\
    &= \displaystyle\sum_{n=0}^{\infty}v_n' z^{2n+1}+\displaystyle\sum_{n=0}^{\infty}v_n'' z^{2n+1}.
\end{aligned}
\end{equation}
while $v_n'$ and $v_n''$ could be directly calculated through Corollary \ref{col2.3}:
\begin{equation}\label{equ3.34}
    \begin{aligned}
        &v_0' = v_0'' = 1, v_1' = \frac{5}{6}, v_1'' = \frac{7}{6}\\
        &v_{n+1}' = \frac{2n^2+3n+\frac{5}{4}}{(n+1)(n+\frac{3}{2})}v_n' - \frac{n(2n+1)}{(n+1)(2n+3)}v_{n-1}'\\
        &v_{n+1}'' = \frac{2n^2+\frac{7}{2}n+\frac{7}{4}}{(n+1)(n+\frac{3}{2})}v_n'' - \frac{(n+\frac{1}{2})^{2}}{(n+1)(n+\frac{3}{2})}v_{n-1}''
    \end{aligned}
\end{equation}

Now, by the Maclaurin expansion of above two special functions, we know that:
\begin{equation*}
    \begin{aligned}
        & arcsinz \sim z + \frac{z^3}{6} + \frac{3}{40}z^5 + o(z^5)\\
        & arctanhz \sim z +\frac{z^3}{3} + \frac{1}{5}z^5 + o(z^5)
    \end{aligned}
\end{equation*}
Therefore, we can just write the $f(z)$ as a power series and differentiate it:
\begin{equation}\label{equ3.35}
    \begin{aligned}
        & f(z) = \displaystyle\sum_{n=0}^{\infty}v_n z^{2n+2}\\
        & f'(z) = \displaystyle\sum_{n=0}^{\infty}(2n+2)v_n z^{2n+1}
    \end{aligned}
\end{equation}
Now compare the coefficients of the $z^{2n+1}$ in (\ref{equ3.33}) and
(\ref{equ3.35}), we have:
\begin{equation}\label{equ3.37}
    v_n = \frac{v_n'+v_n''}{2n+2}
\end{equation}
Combine the results of (\ref{equ3.34}) and (\ref{equ3.37}), we get the final results, which completes the proof.

\end{proof}
\end{corollary}

%\rightarrow
%\int_{0}^{\infty}dt
%\begin{center}

%\end{center}

%\displaystyle\sum_{n=0}^{\infty}
%\prod_{i=1}^{n}
%\genfrac{[}{]}{0pt}{}{2n}{n-1}_q

%
%\bibliographystyle{plain}  % 或者 plainnat、unsrt 等
%\bibliography{sample}

\end{document}